\documentclass[a4paper,fleqn]{cas-sc}

\usepackage{algorithm}
\usepackage{algpseudocode}

\ExplSyntaxOn
\cs_gset:Npn \__first_footerline:
  { \group_begin: \small \sffamily \__short_authors: \group_end: }
\ExplSyntaxOff 
\usepackage{xcolor}
\usepackage{soul}

\usepackage{color,soul}

\usepackage{caption}
\usepackage{booktabs}
\usepackage{threeparttable}

\usepackage{subcaption}

\usepackage{graphicx}

\usepackage{multirow}
\usepackage{amssymb}
\usepackage{amsthm}
\newtheorem{theorem}{Theorem}[section]
\newtheorem{lemma}[theorem]{Lemma}

\usepackage{hyperref}
\usepackage{amsmath}
\usepackage{enumerate}
 \usepackage{graphicx}

\newcommand{\cs}[1]{\texttt{\symbol{`\\}#1}}

\newcounter{example}[section]
\newenvironment{example}[1][]{\refstepcounter{example}\par\medskip
   \noindent \textbf{Example~\theexample. #1} \rmfamily}{\medskip}

\usepackage[numbers]{natbib}
\def\tsc#1{\csdef{#1}{\textsc{\lowercase{#1}}\xspace}}
\tsc{WGM}
\tsc{QE}
\makeatletter

\def\@evenfoot{\normalfont\hfil\thepage\hfil}
\def\@oddfoot{\normalfont\hfil\thepage\hfil}
\makeatother

\usepackage[figcolor=white]{todonotes}

 \usepackage{xcolor}
\begin{document}
%\let\WriteBookmarks\relax
%\def\floatpagepagefraction{1}
%\def\textpagefraction{.001}
% remove comment for short-title

\shorttitle{Structure-preserving mass-lumped FEM for the Allen-Cahn equation}

\shortauthors{H. K. Dwivedi and M. Ehrhardt}
%\begin{frontmatter}

 \title [mode = title]{Mass-lumped finite elements with maximum-bound preservation and original-energy dissipation for the degenerate-mobility Allen-Cahn equation}
 
%\tnotemark[1,2]
%
%\tnotetext[1]{This document is the results of the research project funded by .....}
%
%\tnotetext[2]{The second title footnote which is a longer text matter to fill through the whole text width and overflow into another line in the footnotes area of the first page.}

\author[BHU]{Himanshu Kumar Dwivedi}[
                        orcid=0000-0002-2829-8764]
\ead{himanshukrdwivedi.rs.mat21@itbhu.ac.in}
%\cormark[1]

\address[BHU]{Department of Mathematical Sciences, Indian Institute of Technology (BHU), Varanasi-221005, India.}
%\cormark[a]

\author[BUW]{Matthias Ehrhardt}[
                   orcid=  0000-0003-2561-8854]         
\ead{ehrhardt@uni-wuppertal.de}
\cormark[1]
\cortext[cor1]{Corresponding author}
\address[BUW]{Applied and Computational Mathematics, University of Wuppertal, Gau\ss stra\ss e 20, 42119 Wuppertal, Germany }

%%%%%%%%%%%%%%%%%%%%%%%%%%% abstract
\begin{abstract}
In this work, we propose a mass-lumped finite element method for the Allen-Cahn equation that incorporates nonlinear mobility which may degenerate. 
The spatial approximation uses continuous, piecewise linear, mass-lumped finite elements, while convex splitting is employed for the double-well potential in time. 
Under the standard non-obtuse mesh condition, it is shown that the resulting fully discrete method preserves the maximum-bound principle unconditionally: if the initial approximation takes values in $[-1, 1]$, then so does the numerical solution at every time level. 
Additionally, the method dissipates the original discrete Ginzburg-Landau energy without any time-step restriction, even when the mobility vanishes at pure phases. 
This analysis relies on the monotonicity of the lumped finite element operator and the convexity of the split potential. 
Numerical experiments confirm these theoretical properties and demonstrate the robustness of the method in phase-separation dynamics involving degenerate mobility.
%
%In this work, we propose a mass-lumped finite element method for the Allen-Cahn equation incorporating a nonlinear mobility that may degenerate. 
%The spatial approximation is based on continuous piecewise linear finite elements with mass lumping, while a convex-splitting treatment of the double-well potential is employed in time. 
%Under the standard non-obtuse mesh condition, the resulting fully discrete method is shown to preserve the maximum-bound principle unconditionally:  if the initial approximation takes values in $[-1,1]$, then so does the numerical solution at every time level. 
%In addition, the method dissipates the original discrete Ginzburg-Landau energy without any time-step restriction, even when the mobility vanishes at the pure phases. 
%The analysis relies on the monotonicity structure of the lumped finite element operator and the convexity of the split potential. 
% Numerical experiments confirm the theoretical properties and illustrate the robustness of the method for phase-separation dynamics with degenerate mobility.
\end{abstract}

\begin{keywords}
Allen-Cahn equation \sep degenerate mobility \sep mass-lumped finite element method \sep maximum-bound principle \sep energy dissipation
\end{keywords}

\maketitle

%%%%%%%%%%%%%%%%%%%%%%%%%%%%%%%%%%%%%%%%% Intro
\section{Introduction}
The Allen-Cahn equation is a fundamental phase-field model used to describe interfacial dynamics and phase separation. It was originally introduced to describe antiphase boundary motion in crystalline solids \cite{AllenCahn1979}.
A classic numerical analysis of the model and its relation to mean-curvature flow can be found in \cite{FengProhl2003}. 
In many applications, the relaxation rate depends on the local phase state and is therefore modeled through as a nonlinear mobility. 
In this paper, we investigate the \textit{generalized Allen-Cahn model} described by the following system:
\begin{equation}\label{eq:model}
\begin{cases}
\frac{\partial u}{\partial t} = -M(u)\mu, \qquad
\mu = -\varepsilon^2\Delta u + f(u),
& (\mathbf{x},t)\in\Omega\times(0,T], \\
\frac{\partial u}{\partial\mathbf{n}} = 0,
& (\mathbf{x},t)\in\partial\Omega\times(0,T], \\
u(\mathbf{x},0) = u_0(\mathbf{x}),
& \mathbf{x}\in\Omega.
\end{cases}
\end{equation}
Here, $\Omega\subset\mathbb{R}^d$ $(d=2,3)$ is a bounded polygonal or polyhedral domain, $\mathbf n$ represents the outward-pointing unit normal vector, and $\varepsilon>0$ is the interfacial width. Hereafter, $(\cdot,\cdot)$ and $\|\cdot\|_0$ denote the $L^2(\Omega)$ inner product and norm, respectively. Moreover,
\begin{equation*}
  F(u) = \frac{(u^2-1)^2}{4},\quad \text{is the bistable potential with }\; f(u)=F'(u)=u^3-u.
\end{equation*}
The mobility $M(u)$ is assumed to be nonnegative for $u\in[-1,1]$ and is allowed to degenerate at the pure phases. 
One example is $M(u)=(1-u^2)^m$, where $m$ is a positive integer. 
The considered model has two important structural properties. 
If $\|u_0\|_{L^\infty(\Omega)}\le1$, the solution remains bounded by the pure phases, while the \textit{Ginzburg-Landau energy}
\begin{equation}\label{eq_energy}
  E(u) =\int_\Omega \Bigl(\frac{\varepsilon^2}{2}|\nabla u|^2+F(u)\Bigr)\,d\mathbf{x}
\end{equation}
formally satisfies
\begin{equation*}
     \frac{dE}{dt} = -\bigl(M(u)\mu,\mu\bigr)\le0.
\end{equation*}
%The necessity of maintaining both properties at the discrete level is even more crucial when $M$ is degenerate, as the violation of the maximum-bound can easily break down the nonnegativity and physical interpretation of the mobility.
Maintaining both properties at the discrete level is even more crucial when $M$ is degenerate because violating the maximum-bound can easily cause the mobility to become negative and lose its physical interpretation.

Considerable work has been devoted to structure-preserving discretizations of Allen-Cahn-type equations. 
For nonlinear and possibly degenerate mobility, Shen, Tang and Yang \cite{ShenTangYang2016} constructed maximum-principle-preserving finite-difference (FD) schemes for the generalized Allen-Cahn equation.
Meanwhile, Du and Hou \cite{DuHou2024} developed a second-order, maximum-bound-preserving FD %finite-difference 
method for a generalized model with nonlinear mobility and convection. 
Second-order BDF and stabilized Crank-Nicolson schemes for general mobility were subsequently analyzed in \cite{HouJuQiao2023}. 
%HouJuQiao2024}. 
More recently, Hou, Liu and Ju \cite{HouLiuJu2025} obtained a second-order 
% finite-difference 
FD method that simultaneously preserves the maximum-bound principle and dissipates the original energy. 
Hou, Zhang and Zhu \cite{HouZhangZhu2025} developed an unconditionally maximum-bound-preserving second-order scheme. 
Fu et al.\ \cite{FuHouQiaoZhang2026}  addressed degenerate mobilities through dynamic stabilization. 
Their first-order method unconditionally preserves both the maximum-bound principle and energy dissipation, while the second-order formulation satisfies a conditional maximum-bound principle and a uniform original energy bound. 

The finite element method (FEM) provides greater geometric flexibility, but preserving these structures simultaneously is more challenging. 
Barrett and Blowey \cite{BarrettBlowey2001} developed an early finite element analysis for degenerate Allen-Cahn/Cahn-Hilliard systems. 
For the classical Allen-Cahn equation, energy-stable FEMs have been studied in \cite{LiHuangYi2019}. Mass-lumped, maximum-principle-preserving formulations have also been developed for surface Allen-Cahn equations \cite{XiaoHeFeng2020}. 
%%%
Subsequently, Yang, Yi, and Zhang \cite{YangYiZhang2023} proposed high-order mass-lumped FEM schemes that unconditionally preserve the maximum principle for the classical Allen-Cahn equation. 
However, FEM theory for nonlinear mobility, which can vanish at pure phases, is comparatively underdeveloped. 
This is especially true when a rigorous nodal maximum-bound and dissipation of the original Ginzburg-Landau energy are both required within the same fully discrete formulation.

Motivated by this observation, we construct a mass-lumped $P_1$ finite element scheme for \eqref{eq:model} in conjunction with convex-splitting for the double-well potential.
 The diagonal property of mass lumping, combined with the sign structure of the stiffness matrix on non-obtuse simplicial meshes,  enables us to establish a maximum-bound principle for the scheme using a direct nodal argument.  
Meanwhile, the convex splitting results in the original discrete energy law. 
We prove that, for arbitrary time step sizes, the method satisfies
\begin{equation*}
    |U_h^n|_{L^\infty(\Omega)}\le1,
     \quad\text{and}\quad
E_h(U_h^n)\le E_h(U_h^{n-1}),
 \quad\text{provided}\quad \|U_h^0\|_{L^\infty(\Omega)}\le1.
\end{equation*}
%provided $\|U_h^0\|_{L^\infty(\Omega)}\le1$. 
In particular, 
the analysis remains valid when $M(\pm1)=0$. 
We conducted numerical experiments to assess the accuracy of the proposed scheme, confirming that it maintains both the maximum-bound principle and the energy-dissipation property.

%%%%%%%%%%%%%%%%%%%%% Section 2
\section{Mass-lumped finite element approximation} \label{sec:scheme}
Let $\mathcal{T}_h$ be a conforming and shape-regular simplicial partition of $\Omega$, consisting of triangles for $d=2$ and tetrahedra for $d=3$, and let $h=\max_{K\in\mathcal{T}_h}\operatorname{diam}(K)$. 
We consider the continuous piecewise linear finite element space
\begin{equation*}
   V_h = \bigl\{v_h\in C(\overline{\Omega})\colon\,v_h|_K\in\mathbb{P}_1(K)
\quad\forall\,K\in\mathcal{T}_h \bigr\}.
\end{equation*}
%We denote the vertices of the triangulation $\mathcal{T}_h$ by 
%% the set of indices 
%$\{x_i\}_{i=1}^{N_h}$, and the corresponding nodal basis 
%satisfying $\phi_i(x_j)=\delta_{ij}$.
We denote by $\{x_i\}_{i=1}^{N_h}$ the vertices of $\mathcal{T}_h$ and by $\{\phi_i\}_{i=1}^{N_h}$ the corresponding nodal basis satisfying $\phi_i(x_j)=\delta_{ij}$.
For any nodal functions $v$ and $w$, we introduce the \textit{mass-lumped inner product}
\begin{equation}\label{eq:lumped_product}
   (v,w)_h := \sum_{i=1}^{N_h} v(x_i)\,w(x_i)\,m_i, \quad
   \text{where}\quad m_i := \int_\Omega\phi_i(\mathbf{x})\,d\mathbf{x}>0,
\end{equation}
and introduce the induced norm $\|v_h\|_h^2=(v_h,v_h)_h$.
Thus, the associated mass matrix is diagonal with strictly positive diagonal entries. We recall that we also use the standard inner product (without the subscript  $h$).
Mass lumping is particularly useful in maximum-principle analysis since it preserves the nodal character of the discrete evolution \cite{LiHuang2013, YangYiZhang2023}.

%%%%%%%%%%%%%% stiffness matrix
For the stiffness matrix $A=(a_{ij})_{i,j=1}^{N_h}$, let $a_{ij} := (\nabla\phi_j,\nabla\phi_i)$. 
Throughout the analysis, we impose the standard non-obtuse mesh condition
\begin{equation}\label{eq:mesh_condition}
     a_{ij}\le0,\qquad i\neq j.
\end{equation}
This condition \eqref{eq:mesh_condition} ensures that the finite element stiffness matrix $A$ possesses the monotonicity structure required by a discrete maximum principle % \cite{CiarletRaviart1973}
\cite{LiHuang2013}. 
Since the nodal basis forms a partition of unity and no essential boundary degrees of freedom are eliminated under the homogeneous Neumann condition, we have
\begin{equation}\label{eq:row_sum}
    \sum_{j=1}^{N_h} a_{ij} = 0, \qquad 1\le i\le N_h.
\end{equation}

%%%% energy-stable temporal approximation
To obtain an \textit{energy-stable temporal approximation}, we use a convex splitting of the double-well potential
\begin{equation}\label{eq:convex_split}
    F(s) = F_c(s) - F_e(s), \quad \text{with}\qquad
     F_c(s) = \tfrac{1}{4}s^4 + \tfrac{1}{4},\qquad
    F_e(s) = \tfrac{1}{2}s^2,
\end{equation}
where both $F_c$ and $F_e$ are convex. 
Consequently, $F_c'(s)=s^3,$ and $F_e'(s)=s$. 
This decomposition enables us to treat the convex contribution implicitly and the concave contribution explicitly. This is the key mechanism behind the discrete energy law \cite{Eyre1998, LiHuangYi2019}.

Next, let $0=t_0<t_1<\cdots<t_N = T$ be a uniform temporal partition with time step $\tau=T/N$, and define
\begin{equation*}
    \delta_\tau U_h^{n+1} := \frac{U_h^{n+1}-U_h^n}{\tau}.
\end{equation*}
Assuming $u_0\in C(\overline\Omega)$, we set $U_h^0=I_hu_0$, where $I_h$ denotes the nodal interpolation operator. 
Given $U_h^n\in V_h$, the proposed mass-lumped FEM seeks to find $(U_h^{n+1},\mu_h^{n+1})\in V_h\times V_h$ such that, for all $v_h,q_h\in V_h$,
\begin{subequations}\label{eq:scheme}
\begin{align}
\bigl( \delta_\tau U_h^{n+1}, v_h \bigr)_h +
\bigl( M(U_h^n)\mu_h^{n+1}, v_h \bigr)_h &=0, \label{eq:scheme1}\\
\bigl( \mu_h^{n+1}, q_h \bigr)_h
&=\varepsilon^2 \bigl( \nabla U_h^{n+1}, \nabla q_h \bigr)
+ \bigl( (U_h^{n+1})^3-U_h^n, q_h \bigr)_h.\label{eq:scheme2}
\end{align}
\end{subequations}
Here and below, nonlinear quantities appearing in $(\cdot,\cdot)_h$ are evaluated nodally. 
In particular, $M_i^n := M\bigl(U_h^n(x_i)\bigr)$.
No strictly positive lower bound on the mobility is imposed. We only require
\begin{equation}\label{eq:mobility_assumption}
    M(s)\ge0,\qquad s\in[-1,1],
\end{equation}
so that genuinely degenerate mobilities, such as
\begin{equation*}
     M(s)=(1-s^2)^m,\qquad m\in\mathbb{N},
\end{equation*}
are included. 
Notice that if $M_i^n=0$, then equation \eqref{eq:scheme1} simply implies $U_h^{n+1}(x_i) = U_h^n(x_i)$ at the corresponding node. Therefore, the discrete formulation remains well-defined at  pure phases.

%%%%%%%%%%%%%%%%%%%%%%%%%%%%%%%
For later use, we record the nodal representation of \eqref{eq:scheme}.
%\eqref{eq:scheme1}--\eqref{eq:scheme2}.
Choosing $v_h=q_h=\phi_i$ yields the following system: 
\begin{subequations} \label{eq:nodal}
\begin{align}
  m_i\frac{U_i^{n+1}-U_i^n}{\tau} + m_i M_i^n \mu_i^{n+1} &=0, \label{eq:nodal1}\\
  m_i \mu_i^{n+1} &= \varepsilon^2\sum_{j=1}^{N_h}a_{ij}U_j^{n+1} + m_i\big[(U_i^{n+1})^3-U_i^n\big], \label{eq:nodal2}
\end{align}
\end{subequations}
where $U_i^n:=U_h^n(x_i),$ and $\mu_i^n:=\mu_h^n(x_i).$ Eliminating $\mu_i^{n+1}$ from \eqref{eq:nodal}, 
%\eqref{eq:nodal1}-\eqref{eq:nodal2}, 
we obtain
\begin{equation}\label{eq:nodal_scheme}
   U_i^{n+1} + \tau M_i^n \biggl[ \frac{\varepsilon^2}{m_i} \sum_{j=1}^{N_h}a_{ij}\,U_j^{n+1} + (U_i^{n+1})^3 \biggr] 
   = (1+\tau M_i^n)\,U_i^n .
\end{equation}
Together with \eqref{eq:mesh_condition}--\eqref{eq:row_sum}, the representation \eqref{eq:nodal_scheme} will provide a direct nodal mechanism for establishing the maximum-bound principle.
%
%The representation \eqref{eq:nodal_scheme}, together with \eqref{eq:mesh_condition}--\eqref{eq:row_sum}, will provide a direct nodal mechanism for establishing the maximum-bound principle.

%%%%%%%%%%%%%%%%%%%%%%%%% mass-lumped discrete energy
Finally, corresponding to the Ginzburg-Landau energy in \eqref{eq_energy}, we define the \textit{mass-lumped discrete energy}
\begin{equation}\label{eq:discrete_energy}
   E_h(V_h) := \frac{\varepsilon^2}{2}\|\nabla V_h\|_0^2 + \bigl(F(V_h),1\bigr)_h
   = \frac{\varepsilon^2}{2}\|\nabla V_h\|_0^2 + \sum_{i=1}^{N_h} m_i F(V_h(x_i)).
\end{equation}
The combination of the diagonal lumped mass matrix and the stiffness-matrix sign structure yields the discrete maximum-bound. 
Meanwhile, the convex splitting in \eqref{eq:scheme2} leads directly to dissipation of the original energy \eqref{eq:discrete_energy}. 
These two essential properties are established in the next section.

%%%%%%%%%%%%%%%%%%%%%%%%%%%%%%%%%%%%%%%%%%%%%%%% Section 3
\section{Maximum-bound and energy-dissipation properties} \label{sec:analysis}
We will now establish the two main structural properties of the scheme \eqref{eq:scheme}.
%\eqref{eq:scheme1}--\eqref{eq:scheme2}. 
Throughout this section, we will assume that the mesh condition~\eqref{eq:mesh_condition} and the mobility assumption~\eqref{eq:mobility_assumption} hold.

% \medskip\noindent\textbf{Well-posedness.}
\subsection{Well-posedness}
Suppose that $M_i^n:=M(U_i^n)\ge0$ for all $i$. If $M_i^n=0$, then \eqref{eq:nodal1} immediately yields  $U_i^{n+1}=U_i^n$. 
On the remaining set $\mathcal{I}_+^n:=\{i:M_i^n>0\}$, eliminating $\mu_i^{n+1}$ from
\eqref{eq:nodal}
% \eqref{eq:nodal1}--\eqref{eq:nodal2} 
gives
\begin{equation}\label{eq:active_system}
\frac{m_i}{\tau M_i^n}(U_i^{n+1}-U_i^n)
+\varepsilon^2\sum_{j=1}^{N_h}a_{ij}U_j^{n+1}
+m_i\big[(U_i^{n+1})^3-U_i^n\big]=0 .
\end{equation}
The system \eqref{eq:active_system} is the Euler equation of the functional
\begin{equation}\label{eq:functional}
\mathcal{J}_n(V_h) = \sum_{i\in\mathcal{I}_+^n}
\frac{m_i}{2\tau M_i^n}(V_i-U_i^n)^2
+\frac{\varepsilon^2}{2}\|\nabla V_h\|_0^2
+\sum_{i\in\mathcal{I}_+^n} m_i \Bigl(\frac{1}{4} V_i^4 - U_i^n V_i \Bigr),
\end{equation}
defined on the affine subspace satisfying $V_i=U_i^n$ for $i\notin\mathcal{I}_+^n$. 
Since the first term in \eqref{eq:functional} is strictly convex on the active nodes and the remaining quadratic and quartic terms are convex, $\mathcal{J}_n$ is coercive and strictly convex. 
Therefore, it admits a unique minimizer. 
The chemical potential is subsequently uniquely determined  from \eqref{eq:nodal2}.
Therefore, the fully discrete problem is well-posed.

%%%%%%%%%%%%%%%%%%%%%%%%%% Section 3.2.
\subsection{Maximum-bound analysis}
The following elementary property of the stiffness operator is central to the maximum-bound analysis.
\begin{lemma}\label{lem:extremum}
Let $V_h=\sum_{j=1}^{N_h}V_j\phi_j\in V_h$ and assume
\eqref{eq:mesh_condition}--\eqref{eq:row_sum}. 
If $V_i=\max_jV_j$, then $\sum_{j=1}^{N_h}a_{ij}V_j\ge0.$
%\begin{equation}\label{eq:max_stiffness}
%\sum_{j=1}^{N_h}a_{ij}V_j\ge0.
%\end{equation}
If $V_i=\min_jV_j$, then $\sum_{j=1}^{N_h}a_{ij}V_j\le0.$
%\begin{equation}\label{eq:min_stiffness}
%\sum_{j=1}^{N_h}a_{ij}V_j\le0.
%\end{equation}
\end{lemma}
%%%%%%%%%%%%%%%%%
\begin{proof}
Using \eqref{eq:row_sum}, we have
\begin{equation*}
\sum_{j=1}^{N_h}a_{ij}V_j
=
\sum_{j\neq i}a_{ij}(V_j-V_i).
\end{equation*}
Since $a_{ij}\le0$ for $j\neq i$, and at a nodal maximum, $V_j-V_i\le0$, the first assertion is proven. 
The argument at a nodal minimum is identical, but with the inequalities reversed.
\end{proof}
%%%%%%%%%%%%%%%%%%%%%%

We are now ready to establish the \textit{discrete maximum-bound principle}.
\begin{theorem}[maximum-bound principle] \label{thm:mbp}
Assume that $\|U_h^0\|_{L^\infty(\Omega)}\le1,$ and that $M(s)\ge0$ for $s\in[-1,1]$. 
Then, for arbitrary $\tau>0$, the solution of \eqref{eq:scheme}
%\eqref{eq:scheme1}--\eqref{eq:scheme2} 
satisfies
\begin{equation}\label{eq:MBP}
    \|U_h^n\|_{L^\infty(\Omega)}\le1, \qquad 0\le n\le N.
\end{equation}
\end{theorem}
%%%%%%%%%%%%%%%%%%%
\begin{proof}
We proceed by induction. 
The assertion is true for $n=0$. Suppose that $|U_i^n|\le1$ for all nodes $x_i$. 
By \eqref{eq:mobility_assumption}, we have $M_i^n=M(U_i^n)\ge0.$ Assume, to the contrary, that
\begin{equation*}
   U_{i_*}^{n+1} = \max_{1\le i\le N_h}U_i^{n+1}>1.
\end{equation*}
By Lemma~\ref{lem:extremum},
\begin{equation*}
    \sum_{j=1}^{N_h} a_{i_*j} \,U_j^{n+1}\ge0.
\end{equation*}
Hence, from the nodal equation \eqref{eq:nodal_scheme},
\begin{align*}
(1+\tau M_{i_*}^n)\,U_{i_*}^n &= U_{i_*}^{n+1} +\tau M_{i_*}^n
\biggl[ \frac{\varepsilon^2}{m_{i_*}}
\sum_{j=1}^{N_h}a_{i_*j}\,U_j^{n+1} +(U_{i_*}^{n+1})^3\biggr]  \ge U_{i_*}^{n+1} 
  + \tau M_{i_*}^n\,(U_{i_*}^{n+1})^3 > 1+\tau M_{i_*}^n . 
\end{align*}
On the other hand, the induction hypothesis gives
\begin{equation*}
    (1+\tau M_{i_*}^n)\,U_{i_*}^n  \le 1+\tau M_{i_*}^n,
\end{equation*}
which is a contradiction. Therefore, $\max_i U_i^{n+1}\le1$. 
Similarly, suppose that
\begin{equation*}
    U_{i_*}^{n+1} = \min_{1\le i\le N_h} U_i^{n+1} < -1.
\end{equation*}
From Lemma~\ref{lem:extremum}  we get the following
\begin{equation*}
   \sum_{j=1}^{N_h} a_{i_*j}\,U_j^{n+1}\le0.
\end{equation*}
Using \eqref{eq:nodal_scheme}, we obtain
\begin{align*}
  (1+\tau M_{i_*}^n)\,U_{i_*}^n
  &\le U_{i_*}^{n+1} +\tau M_{i_*}^n(U_{i_*}^{n+1})^3
   < -1-\tau M_{i_*}^n,
\end{align*}
whereas the induction hypothesis yields
\begin{equation*}
    (1+\tau M_{i_*}^n)\,U_{i_*}^n \ge -1-\tau M_{i_*}^n.
\end{equation*}
This is again a contradiction. 
Thus $-1\le U_i^{n+1}\le1$ at every node. 
Since $U_h^{n+1}$ is piecewise linear, its extrema on every simplex are attained at the vertices, and hence \eqref{eq:MBP} follows. 
The induction is complete.
\end{proof}
%%%%%%%%%%%%%%%%%%%%%%%%%

% An immediate 
A consequence of Theorem~\ref{thm:mbp} is that
$M(U_h^n)\ge0$ at every time level, including the degenerate case $M(\pm1)=0$.

%%%%%%%%%%%%%%%%%%%% Section 3.3 %%%%%
\subsection{Energy-dissipation properties}
We next show that the same scheme dissipates the original discrete Ginzburg-Landau energy \eqref{eq:discrete_energy}.
%%%%%%%%%%%%%%%%%%%%%%%%
\begin{theorem}[Original-energy dissipation]\label{thm:energy}
Under the assumptions of Theorem~\ref{thm:mbp}, the solution of \eqref{eq:scheme}
%\eqref{eq:scheme1}--\eqref{eq:scheme2} 
satisfies the discrete energy identity
\begin{equation}\label{eq:energy_identity}
E_h(U_h^{n+1}) - E_h(U_h^n)
  =-\tau \bigl(M(U_h^n)\mu_h^{n+1},\mu_h^{n+1}\bigr)_h
   - \frac{\varepsilon^2}{2}\|\nabla(U_h^{n+1}-U_h^n)\|_0^2
-\mathcal{R}_h^{n+1},
\end{equation}
where
\begin{equation}\label{eq:energy_remainder}
\mathcal{R}_h^{n+1} := \frac{1}{4}\sum_{i=1}^{N_h}m_i
(U_i^{n+1}-U_i^n)^2
\Bigl[ 3(U_i^{n+1})^2 +2U_i^{n+1}U_i^n +(U_i^n)^2+2 \Bigr]\ge0 .
\end{equation}
Consequently, for every $\tau>0$, we have the discrete energy decay
\begin{equation}\label{eq:energy_decay}
    E_h(U_h^{n+1}) \le E_h(U_h^n),
\qquad 0\le n\le N-1.
\end{equation}
%for every $\tau>0$.
\end{theorem}
%%%%%%%%%%%%%%%%%%%%
\begin{proof}
First, we set $\delta U_h^{n+1}:=U_h^{n+1}-U_h^n.$ Taking $v_h=\mu_h^{n+1}$ in \eqref{eq:scheme1} and multiplying by $\tau$ gives
\begin{equation}\label{eq:energy_step1}
\bigl( \delta U_h^{n+1}, \mu_h^{n+1} \bigr)_h
+ \tau\bigl( M(U_h^n)\mu_h^{n+1}, \mu_h^{n+1}\bigr)_h = 0.
\end{equation}
Taking $q_h=\delta U_h^{n+1}$ in \eqref{eq:scheme2}, we obtain
\begin{equation}\label{eq:energy_step2}
\bigl( \mu_h^{n+1},\delta U_h^{n+1} \bigr)_h
= \varepsilon^2 \bigl( \nabla U_h^{n+1}, \nabla\delta U_h^{n+1} \bigr) +
\bigl((U_h^{n+1})^3-U_h^n, \delta U_h^{n+1}\bigr)_h .
\end{equation}
The gradient contribution satisfies the exact identity
\begin{equation}\label{eq:gradient_identity}
  \bigl( \nabla U_h^{n+1}, \nabla\delta U_h^{n+1} \bigr)
  = \frac{1}{2} \bigl( \|\nabla U_h^{n+1}\|_0^2 -\|\nabla U_h^n\|_0^2
   +\|\nabla\delta U_h^{n+1}\|_0^2 \bigr).
\end{equation}
Moreover, for arbitrary $a,b\in\mathbb R$ and $F(s)=\frac{1}{4}(s^2-1)^2$, it is straightforward that
\begin{equation}\label{eq:potential_identity}
   (a^3-b)(a-b) = F(a)-F(b) + \frac{(a-b)^2}{4} \bigl(3a^2+2ab+b^2+2\bigr).
\end{equation}
Since
\begin{equation*}
   3a^2+2ab+b^2+2 = 2a^2+(a+b)^2+2 > 0,
\end{equation*}
the last term in \eqref{eq:potential_identity} is nonnegative. 
Applying \eqref{eq:potential_identity} nodewise in the lumped inner product gives
\begin{equation}\label{eq:potential_discrete}
\bigl((U_h^{n+1})^3-U_h^n, \delta U_h^{n+1}\bigr)_h = \bigl(F(U_h^{n+1})-F(U_h^n),1\bigr)_h +\mathcal{R}_h^{n+1}.
\end{equation}
Substitution of \eqref{eq:gradient_identity} and
\eqref{eq:potential_discrete} into \eqref{eq:energy_step2}, followed by \eqref{eq:energy_step1}, yields \eqref{eq:energy_identity}. 
Finally, Theorem~\ref{thm:mbp} ensures $M(U_i^n)\ge0$ for every node, and therefore
\begin{equation*}
\bigl(M(U_h^n)\mu_h^{n+1},\mu_h^{n+1}\bigr)_h
= \sum_{i=1}^{N_h} m_iM(U_i^n)(\mu_i^{n+1})^2 \ge0.
\end{equation*}
Then, the claimed energy dissipation property follows, along with \eqref{eq:energy_remainder}, which implies \eqref{eq:energy_decay}.
\end{proof}
%%%%%%%%%%%%%%%%%%%
Thus, the proposed mass-lumped FEM preserves the intrinsic structures of the degenerate-mobility Allen-Cahn dynamics without any restriction on the time-step size. 
In particular, the maximum-bound keeps the numerical solution within the physically admissible interval $[-1,1]$, while \eqref{eq:energy_identity} shows that the original discrete Ginzburg-Landau energy, rather than a modified one, is dissipated at every time step.

%%%%%%%%%%%%%%%%%% Section 4 %%% Numerics
\section{Numerical experiments}\label{sec:numerics}
In this section, we present two numerical examples that demonstrate the accuracy and structure-preserving properties of the proposed method. 
The computational domain is triangulated by dividing each uniform rectangular cell into two right triangles.
Thus, the resulting mesh satisfies the non-obtuse condition \eqref{eq:mesh_condition}. 
All nonlinear systems are solved by Newton iteration using $U_h^n$ as the initial guess.
The iterative procedure stops once the relative residual decreases below the prescribed tolerance of $10^{-10}$. 
Unless otherwise stated, we use the degenerate mobility $M(u)=1-u^2$.

%%%%%%%%%%%%%%%%%% Example 1
\begin{example}[(Accuracy test)]\label{ex41}
%\medskip
%\noindent\textbf{Example 4.1 (Accuracy test).}
First, we consider the forced problem
\begin{equation}\label{eq:forced_problem}
   u_t = -M(u)\mu + g(\mathbf{x},t),\qquad \text{with}
  \qquad
   \mu = -\varepsilon^2\Delta u+u^3-u,
\end{equation}
on $\Omega=(0,1)^2$, supplemented with homogeneous Neumann boundary conditions. The source term is chosen such that
\begin{equation}\label{eq:exact_solution}
u(\mathbf{x},t) = \frac{\mathrm{e}^{-t} \cos(\pi x)\,\cos(\pi y)}{2}
\end{equation}
is the exact solution. 
In particular, since $\Delta u=-2\pi^2u,$ we take
\begin{equation}\label{eq:source}
  g = -u + (1-u^2)\, (2\varepsilon^2\pi^2u +u^3 -u).
\end{equation}
The exact solution satisfies the homogeneous Neumann condition and remains strictly within the admissible interval $[-1,1]$.
We set $\varepsilon=0.1$, $T=1$, and choose $\tau=h^2$ so that the first-order temporal error does not dominate the spatial discretization error. At the final time, we measure
\begin{equation*}
  e_0=\|u(\cdot,T)-U_h^N\|_0, \qquad
   e_1=\|\nabla(u(\cdot,T)-U_h^N)\|_0.
\end{equation*}
The corresponding convergence rates are computed by
\begin{equation*}
\mathrm{Rate} = \log_2\Bigl(\frac{e(h)}{e(h/2)}\Bigr).
\end{equation*}
The numerical errors are reported in Table~\ref{tab:convergence}. 
as expected,  the results exhibit the expected second-order convergence in the $L^2$ norm and first-order convergence in the $H^1$ seminorm for the piecewise linear finite element approximation.

%%%%%%%%%%%%%%%%%%%%% Table
\begin{table}[htbp]
\centering
\caption{Numerical errors and convergence rates for Example~\ref{ex41} with $\tau = h^2$.}
\label{tab:convergence}
\renewcommand{\arraystretch}{1.15}
\setlength{\tabcolsep}{7pt}
\begin{tabular}{cccccc}
\toprule
& & \multicolumn{2}{c}{$L^2$-error} 
& \multicolumn{2}{c}{$H^1$-seminorm error} \\
\cmidrule(lr){3-4}\cmidrule(lr){5-6}
$h$ & $\tau$ 
& $\|u(T)-U_h^N\|_{0}$ & Rate
& $\|\nabla(u(T)-U_h^N)\|_{0}$ & Rate \\
\midrule
$1/8$  
& $1.5625\times10^{-2}$ 
& $3.36204\times10^{-3}$ 
& - 
& $8.36695\times10^{-2}$ 
& - \\

$1/16$ 
& $3.9062\times10^{-3}$ 
& $8.69831\times10^{-4}$ 
& $1.951$ 
& $4.05445\times10^{-2}$ 
& $1.045$ \\

$1/32$ 
& $9.7656\times10^{-4}$ 
& $2.19438\times10^{-4}$ 
& $1.987$ 
& $2.01082\times10^{-2}$ 
& $1.012$ \\

$1/64$ 
& $2.4414\times10^{-4}$ 
& $5.49916\times10^{-5}$ 
& $1.997$ 
& $1.00347\times10^{-2}$ 
& $1.003$ \\
\bottomrule
\end{tabular}
\end{table}
\end{example}

%%%%%%%%%%%%%%%%%% Example 2
\begin{example}[(Phase separation with degenerate mobility)]\label{ex42}
%\medskip
%\noindent\textbf{Example 4.2 (Phase separation with degenerate mobility).}
Next, we examine the long-time phase evolution of the unforced Allen-Cahn problem~\eqref{eq:model} on $\Omega=(0,1)^2$. 
We set $\varepsilon=0.02$ with $M(u)=1-u^2$ and initialize the nodal values with a small random perturbation $U_h^0(x_i)=0.05\,\xi_i,$ where $\xi_i\sim\mathcal{U}[-1,1].$
%\begin{equation}\label{eq:random_initial}
%U_h^0(x_i)=0.05\,\xi_i,
%\qquad
%\xi_i\sim\mathcal{U}[-1,1].
%\end{equation}
Hence,
\begin{equation*}
\|U_h^0\|_{L^\infty}\le0.05<1.
\end{equation*}
The computation is performed up to time $T=10$ on a uniform triangulation with $h=1/128$ and time step $\tau=10^{-2}$. 
To examine the two structural properties established in
Theorems~\ref{thm:mbp} and \ref{thm:energy}, we monitor the following:
\begin{equation}\label{eq:diagnostics}
   U_{\max}^n:=\max_{1\le i\le N_h}U_i^n,\qquad
    U_{\min}^n:=\min_{1\le i\le N_h}U_i^n,\qquad
    \mathcal{E}^n:=E_h(U_h^n).
\end{equation}
Representative phase-field configurations, along with the histories of $U_{\max}^n$, $U_{\min}^n$, and $\mathcal{E}^n$ are displayed in Figure~\ref{fig:phase}. 
The nodal extrema remain in the admissible interval, $-1\le U_{\min}^n\le U_{\max}^n\le1$, while the original discrete Ginzburg-Landau energy decreases monotonically, $\mathcal{E}^{n+1}\le\mathcal{E}^n$, in agreement with the theoretical results.
Furthermore, the mobility $M(u)=1-u^2$ decreases as $|u|$ approaches one and vanishes at the pure phases. This illustrates the effect of mobility degeneracy on the interfacial evolution.

%%%%%%%%%%%%%%%%%%%%%%% Figure 1
\begin{figure}
    \centering
 %   \includegraphics[width=0.9\textwidth]{fig1.eps}
     % \includegraphics[width=\textwidth]{fig1.eps}
     %\includegraphics[width=0.8\textwidth]{Fig1.pdf}
     % for arxiv version
     \includegraphics[width=0.9\textwidth]{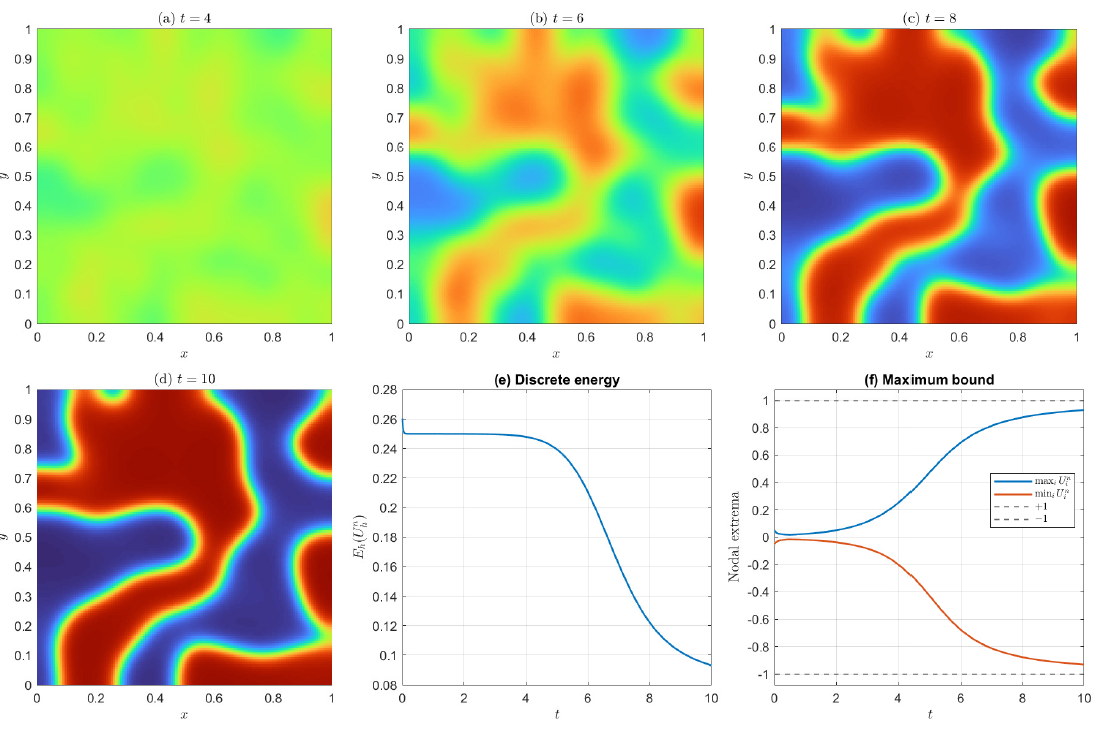}
    \caption{%
    Numerical results for Example~\ref{ex42} with $M(u)=1-u^2$, $\varepsilon=0.02$, $h=1/128$, and $\tau=10^{-2}$.
    Panels (a)--(d) show $U_h^n$ at $t=4,6,8,$ and $10$, respectively. 
    Panel (e) shows the monotone decay of the discrete energy $E_h(U_h^n)$, while panel (f) confirms that the nodal extrema remain in $[-1,1]$, which is in agreement with Theorems~\ref{thm:mbp} and~\ref{thm:energy}.}
    \label{fig:phase}
\end{figure}
\end{example}

%%%%%%%%%%%%%%%%%%%%%%% Conclusion
\section{Conclusion}
We have developed a mass-lumped $P_1$ FEM for the Allen-Cahn problem involving nonlinear, potentially degenerate mobility. 
Combining mass lumping with a convex-splitting treatment of the double-well potential preserves the nodal maximum-bound principle in the proposed scheme, which dissipates the original discrete Ginzburg-Landau energy without any restriction on the time step size. 
The analysis remains valid when the mobility vanishes at the pure phases, and the discrete energy law further reveals the associated numerical dissipation mechanisms. 
Numerical experiments confirmed the expected finite element accuracy and verified the maximum-bound and energy-dissipation properties. 
These results provide a simple, robust, structure-preserving finite element framework for Allen-Cahn dynamics with degenerate mobility.

%%%% Declarations %%%%%%%%%%%%%%%%%%%%
%\section*{Data availability}
%The data used in this work can be obtained from the corresponding author upon justified request.

%\section*{Conflict of Interest}
%The authors affirm that no financial or personal affiliations exist that might have influenced the content or conclusions of this work.

%\section*{Declaration}
%\subsection*{Competing Interest}
%The authors declare that they have no competing interests.

%%%%%%%%%%%%%%%%%%%%%%%%% Refs

\end{document}